\documentclass [a4paper]{IEEEtran}
\usepackage{amsmath}
\usepackage{amsthm}
\usepackage{amssymb, amsfonts}
\usepackage{algorithm}
\usepackage{algorithmic}
\usepackage{graphics}
\usepackage{graphicx}
\usepackage{lipsum}
\usepackage[framed,numbered,autolinebreaks,useliterate,bw]{mcode}
\usepackage[utf8]{inputenc}
\usepackage{skt}

\renewcommand{\u}[1]{\mathrm{#1}}

\title{Coefficient Matrices Computation of Structural Vector Autoregressive Model}
\author{Aravindh Krishnamoorthy\thanks{\hrule}\thanks{\noindent Aravindh Krishnamoorthy is currently with Ericsson Modem Nuremberg GmbH working in the area of Wireless Communication. E-mail: aravindh.krishnamoorthy@ericsson.com, aravindh.k@ieee.org}}

\begin{document}
\maketitle

\begin{abstract}
In this paper we present the Large Inverse Cholesky (LIC) method, an efficient method for computing the coefficient matrices of a Structural Vector Autoregressive (SVAR) model.
\end{abstract}


\section{Introduction}
Structural Vector Autoregressive (SVAR) Model \cite{book:lut} is widely used for multi-branch signal modelling in the fields of Wireless Communication, Econometrics, Physics, and Multi-dimensional Audio Signal Processing.

Given $x = \{x(n)\}$, the input matrix with M signal branches and length N, we would like to express $x$ in the K-th order vector auto-regressive form in terms of time-invariant MxM matrices $L$, $R_i$, and intercept vector $t$, collectively referred to as coefficient matrices, such that the residual $w = \{w(n)\}$ has an identity covariance matrix $I_M$. The model for SVAR is given as follows:
\begin{equation}
	{L} {x}(n) = {t} + \sum_{i=1}^{K} {R}_i {x}(n-i) + {w}(n) \label{eqn:svar}
\end{equation}

In this paper we present an efficient method for computation of the model coefficient matrices. In section \ref{sec:ls} we review a widely used method based on least-squares. In section  \ref{sec:lic} we propose the Large Inverse Cholesky (LIC) method. In section \ref{sec:sc} we compare the computational complexity of these methods, followed by concluding remarks in \ref{sec:con}.

\section{Least Squares Based Method}
\label{sec:ls}
Least-squares based method uses the reduced-form VAR (RVAR) as an intermediate step in computing the coefficient matrices of SVAR. First, the input matrix $x$ is modelled as an RVAR and the coefficient matrices of RVAR are computed using least-squares.
The RVAR model is as follows:
\begin{equation}
	{x}(n) = {c} + \sum_{i=1}^{K} {A}_i {x}(n-i) + {v}(n) \label{eqn:rvar}
\end{equation}
In the above equation, the covariance matrix of $v = \{v(n)\}$ is a positive-definite symmetric matrix.

Equation (\ref{eqn:rvar}) may be rewritten in the matrix-form as follows:
\begin{equation}
	{X} = {A} {S} + {V} \label{eqn:rmat}
\end{equation}
Where, the matrices ${X} \in \mathbb{C}^{\u{Mx(N-K)}}$, ${A} \in \mathbb{C}^{\u{Mx(MK+1)}}$, ${S} \in \mathbb{C}^{\u{(MK+1)x(N-K)}}$, and ${V} \in \mathbb{C}^{\u{MxN-K}}$ are given as:
\begin{align}
{X} = [& {x}(K+1), {x}(K+2), \ldots, {x}(N)] \\
{A} = [& {c}, {A}_1, {A}_2, \ldots, {A}_K] \\
{S} = [& \mathbf{1}_{\u{1x(N-K)}};  \nonumber \\
		&  {x}(K), {x}(K+1), \ldots, {x}(N-1); \nonumber \\
		& {x}(K-1), {x}(K), \ldots, {x}(N-2); \nonumber \\ 
		& \vdots \nonumber \\
		& {x}(1), {x}(2), \ldots, {x}(N-K)] \label{eqn:sdef} \\
{V} = [& {v}(K+1), {v}(K+2), \ldots, {v}(N)]
\end{align}

In the definition of ${S}$, a comma indicates that the next term is concatenated horizontally, thereby increasing the number of columns, while a semicolon indicates that the next term is concatenated vertically, thereby increasing the number of rows.

The parameter matrix of the RVAR model $A$, and the matrix $V$ may be found using the least-squares method as:
\begin{align}
	A &= {X} {S}^H ({S} {S}^H)^{-1} \\
	V &= X - AS
\end{align}

Next, the RVAR model is converted to its equivalent SVAR model by multiplying the RVAR model equation (\ref{eqn:rvar}) with the inverse-Cholesky based whitening filter of $V$. From comparison of equations (\ref{eqn:svar}) and (\ref{eqn:rvar}) we note that the inverse-Cholesky filter is $L$, i.e. $VV^H = L^{-1} (L^{-1})^H$. Therefore, the coefficient matrices of SVAR may be computed from the coefficient matrices of RVAR by multiplying them with $L$.

\section{Large Inverse Cholesky Method}
\label{sec:lic}
Large Inverse Cholesky method computes the coefficient matrices of SVAR directly using an inverse-Cholesky of $TT^H$, where ${T} \in \mathbb{C}^{\u{(M(K+1)+1)x(N-K)}}$ defined as follows:
\begin{align}
{T} = [& \mathbf{1}_{\u{1x(N-K)}};  \nonumber \\
		& {x}(1), {x}(2), \ldots, {x}(N-K); \nonumber \\
		& {x}(2), {x}(3), \ldots, {x}(N-K+1); \nonumber \\
		& \vdots \nonumber \\
		&  {x}(K), {x}(K+1), \ldots, {x}(N-1); \nonumber \\
		&  {x}(K+1), {x}(K+2), \ldots, {x}(N)] \label{eqn:tdef}
\end{align}

Let $U$ be the lower-triangular inverse-Cholesky of $TT^H$, then we have:
\begin{equation}
U^{-1} (U^{-1})^H  = {T} {T}^H \label{eqn:licm}
\end{equation} 

Let index variables:
\begin{align}
	\alpha &= (MK+2, \ldots, M(K+1)+1) \\
	\beta_i &= ((i-1)M+2, \ldots, iM+1)
\end{align}
then, the coefficient matrices of the SVAR model are given by:
\begin{align}
	{L} & = U(\alpha, \alpha) \label{eqn:lic1}\\
	{R_i} & = -U(\alpha, \beta_{K-i+1}) \label{eqn:lic2}\\
	{t} & = -U(\alpha, 1) \label{eqn:lic3}
\end{align}
With $i = 1\ldots K$.

We can easily verify that the coefficient matrices computed by the least-squares method and the LIC method are equivalent. 
Let ${R} = [{t}, {R}_1, {R}_2, \ldots, {R}_K]$, and let the last M rows of $U$ be the matrix $W$, then we have $W = [-R, L]$ and $T = [S; X]$.
From equation (\ref{eqn:licm}) we have:
\begin{equation}
	\begin{bmatrix}
	-R &  L
	\end{bmatrix}
	\begin{bmatrix}
	SS^H &  SX^H \\
	XS^H & XX^H
	\end{bmatrix}
	\begin{bmatrix}
	-R^H \\
	L^H
	\end{bmatrix}
	= I_M
\end{equation}
Upon simplification:
\begin{equation}
	L (X - L^{-1}RS) (X - L^{-1}RS)^H L^H  = I_M
\end{equation}
which is readily identified as the equivalent of least-squares based method.

A MATLAB compatible algorithm without optimization is given below, the full program may be downloaded from the MATLAB Central File Exchange from \cite{url:lic1}.
\begin{lstlisting}
T = zeros(M*(K+1)+1,N-K) ;
T(1,:) = ones(1,N-K) ;
T(M*K+2:end,:) = X(:,K+1:end) ;
for i=1:K
   T(M*(i-1)+2:M*i+1,:) = X(:,i:i-1+N-K) ;
end

U = inv(chol(T*T', 'lower')) ;
alpha = M*K+2:M*(K+1)+1 ;

L = U(alpha, alpha) ;
R = zeros(M,M,K) ;
for i = 1:K
    R(:,:,i) = -U(alpha, (K-i)*M+2:(K-i+1)*M+1) ;
end
t = -U(alpha, 1) ;
\end{lstlisting}

\section{Computational Complexity}
For computing the complexity, we consider the number of multiplies for each operation. 

\newpage

It is assumed that on modern DSP processors, the Multiply-accumulate (MAC) instruction hides the additions that are performed in each step. 
Further, where the output of matrix multiplies result in symmetric matrices, the number of multiplies is chosen as half of the conventional value.
The computational complexity for Cholesky decomposition and inversion is chosen for an NxN matrix at the most efficient value of $N^3/2$.

\label{sec:sc}
\subsection{Least Squares Method}
\begin{tabular}{|p{0.3\columnwidth}|p{0.5\columnwidth}|}
\hline
Operation	&	Number of Multiplies \\
\hline\hline
${S}{S}^H$ 	&	$(\u{MK+1)^2(N-K)}/2$ \\
$({S}{S}^H)^{-1}$	&	$\u{(MK+1)^3}/2$ \\
${X} {S}^H$	&	$\u{M(MK+1)(N-K)}$ \\
${X} {S}^H({S}{S}^H)^{-1}$	&	$\u{M(MK+1)^2}$ \\
${\hat{V}}$	&	$\u{M(MK+1)(N-K)}$ \\
${\hat{V}} {\hat{V}}^H$	&	$\u{M^2(N-K)}$ \\
${L}$	&	$\u{M^3/2}$ \\
${R_i}$		&	$\u{M^3K}$ \\
${t}$	&	$\u{M^2}$ \\
\hline
\end{tabular}

\subsection{Large Inverse Cholesky Method}
\begin{tabular}{|p{0.3\columnwidth}|p{0.5\columnwidth}|}
\hline
Operation	&	Number of Multiplies \\
\hline\hline
${T}{T}^H$ 	&	$(\u{M(K+1)+1)^2 (N-K)}/2$ \\
${\hat{U}}$	&	$\u{(M(K+1)+1)^3/2}$ \\
\hline
\end{tabular}

\vspace{4mm}
LIC is efficient for $\u{M,K} << \u{N}$. For practical cases, LIC can be upto 30\% efficient than least-squares method for finding the coefficient matrices of an SVAR.

\section{Conclusion}
\label{sec:con}
In this paper we presented the Large Inverse Cholesky method for computing the coefficient matrices of a Structural Autoregressive model which is upto 30\% efficient compared to the conventional least-squares based method.


\end{document}